\documentclass[twoside]{article} 
\usepackage{graphicx}
\date{} 
\title{An extension of an asymptotic result of Tricomi concerning a definite integral}
\author{\sc R. B.\ Paris \\
{\em Division of Computing and Mathematics,} \\
{\em Abertay University, Dundee DD1 1HG, UK}}
\begin{document}
\def\f#1#2{\mbox{${\textstyle \frac{#1}{#2}}$}}
\def\dfrac#1#2{\displaystyle{\frac{#1}{#2}}}
\def\boldal{\mbox{\boldmath $\alpha$}}
\newcommand{\bee}{\begin{equation}}
\newcommand{\ee}{\end{equation}}
\newcommand{\lam}{\lambda}
\newcommand{\ka}{\kappa}
\newcommand{\al}{\alpha}
\newcommand{\ba}{\beta}
\newcommand{\la}{\lambda}
\newcommand{\ga}{\gamma}
\newcommand{\eps}{\epsilon}
\newcommand{\fr}{\frac{1}{2}}
\newcommand{\fs}{\f{1}{2}}
\newcommand{\g}{\Gamma}
\newcommand{\br}{\biggr}
\newcommand{\bl}{\biggl}
\newcommand{\ra}{\rightarrow}
\newcommand{\gtwid}{\raisebox{-.8ex}{\mbox{$\stackrel{\textstyle >}{\sim}$}}}
\newcommand{\ltwid}{\raisebox{-.8ex}{\mbox{$\stackrel{\textstyle <}{\sim}$}}}
\renewcommand{\topfraction}{0.9}
\renewcommand{\bottomfraction}{0.9}
\renewcommand{\textfraction}{0.05}
\newcommand{\mcol}{\multicolumn}
\date{}
\maketitle
\pagestyle{myheadings}
\markboth{\hfill \sc R. B.\ Paris  \hfill}
{\hfill \sc Tricomi's integral\hfill}
\begin{abstract}
We consider the expansion of an integral considered by F.G. Tricomi given by
\[\int_{-\infty}^\infty x e^{-x^2}(\f{1}{2}+\f{1}{2}\mbox{erf}\,x)^{m} dx\]
as $m\to\infty$. The procedure involves a suitable change of variable and the inversion of the complementary error function $\mbox{erfc}\,x$. Numerical results are presented to demonstrate the accuracy of the expansion.

A second part examines an extension of an integral arising in airfoil theory.
\vspace{0.3cm}

\vspace{0.1cm}
 
\end{abstract}

\vspace{0.3cm}

\noindent $\,$\hrulefill $\,$

\vspace{0.3cm}

\begin{center}
{\bf 1.\ Introduction}
\end{center}
\setcounter{section}{1}
\setcounter{equation}{0}
\renewcommand{\theequation}{\arabic{section}.\arabic{equation}}
In a short note published in the {\it Accademia dei Lincei} \cite{T}, F.G. Tricomi obtained the leading asymptotic behaviour of the integral arising in probability theory
\[\int_{-\infty}^\infty x e^{-x^2} \bl(\frac{1+\mbox{erf}\,x}{2}\br)^{\!m} dx\]
as $m\to\infty$, where $\mbox{erf}\,x$ is the error function defined by
\[\mbox{erf}\,x=\frac{2}{\sqrt{\pi}} \int_0^x e^{-t^2}dt.\]
The factor in brackets in the above integral satisfies $0<\fs+\fs\mbox{erf}\,x<1$ for $x\in(-\infty,\infty)$, with the upper and lower bounds approached as $x\to\pm\infty$, respectively.
By a suitable transformation of the integration variable, Tricomi showed that the leading behaviour of this integral is given by
\bee\label{e10}
\hspace{1cm}\frac{\sqrt{\pi \log\,m}}{m}\qquad (m\to\infty).
\ee

In the first part of this paper we extend the result in (\ref{e10}) by obtaining higher-order terms in the large-$m$ expansion of the integrals
\bee\label{e11}
I_{n,m}=\int_{-\infty}^\infty x^n e^{-x^2} \bl(\frac{1+\mbox{erf}\,x}{2}\br)^{\!m} dx\qquad (n=0, 1, 2),
\ee
the case $n=0$ being trivial. The second part considers the evaluation of an extension of an integral arising in airfoil theory.
\vspace{0.6cm}

\begin{center}
{\bf 2.\ Extension of Tricomi's analysis}
\end{center}
\setcounter{section}{2}
\setcounter{equation}{0}
\renewcommand{\theequation}{\arabic{section}.\arabic{equation}}
Following Tricomi \cite{T}, we make the change of variable
\bee\label{e21}
e^{-t}=\frac{1+\mbox{erf}\,x}{2},\qquad \frac{e^{-x^2}}{\sqrt{\pi}}\,\frac{dx}{dt}=-e^{-t}
\ee
to yield 
\bee\label{e22}
I_{n,m}=-\int_0^\infty x^n e^{-x^2} \frac{dx}{dt}\,e^{-mt}\,dt=\sqrt{\pi} \int_0^\infty x^n e^{-st}\,dt,\qquad s:=m+1.
\ee
When $n=0$, we then have the trivial result
\[I_{0,m}=\frac{\sqrt{\pi}}{m+1}.\]

To deal with the cases $n=1$ and $n=2$, we require the inversion of the first relation in (\ref{e21}) written in the form
\bee\label{e22a}
\mbox{erfc}\,x=2-2e^{-t}
\ee
to yield $x\equiv x(t)$, which is finite and continuous for $0<t<\infty$.  From the well-known asymptotic behaviour \cite[p.~164]{DLMF}
\[\mbox{erfc}\,x\sim \frac{e^{-x^2}}{\sqrt{\pi} x},\qquad \mbox{erfc} (-x)\sim 2-\frac{e^{-x^2}}{\sqrt{\pi} x}
\qquad (x\to +\infty),\]
it is seen that $x(t)\sim \sqrt{-\log\,t}$ as $t\to 0$ and $x(t)\sim -\sqrt{t}$ as $t\to+\infty$, with $x(t)=0$ when $t=\log\,2$.

To proceed we split the integral into two parts: $[0,\al]$ and $(\al,\infty)$, where $\al$ is finite and bounded away from zero. Then, with the change of variable $t=u/s$, we have
\bee\label{e23}
I_{n,m}=\frac{\sqrt{\pi}}{s} \int_0^{\al s} x^n e^{-u} du+\frac{\sqrt{\pi}}{s} \int_{\al s}^\infty x^n e^{-u} du
\ee
where $x\equiv x(u/s)$. Since $x(t)>-K\sqrt{t}$ for some suitable positive constant $K$, the second integral is bounded by
\[\frac{1}{s}\int_{\al s}^\infty |x(u/s)|^n e^{-u}du<\frac{K^n}{s^{n/2+1}} \int_{\al s}^\infty u^{n/2}e^{-u}du
=\frac{K^n}{s^{n/2+1}} \g(\fs n+1,\al s)=O(s^{-1} e^{-\al s}),\]
as $s\to\infty$, where we have employed the result for the (upper) incomplete gamma function $\g(a,z)\sim z^{a-1} e^{-z}$ as $z\to+\infty$. 

In the first integral in (\ref{e23}) the dominant contribution will arise from the neighbourhood of  $u=0$ as $m\to\infty$. 
We therefore require the inversion of (\ref{e22a}) as $t\to0$ expressed in the form
\bee\label{e23inv}
\mbox{erfc}\,x=\frac{2u}{s}\bl(1-\frac{u}{2s}+\cdots\br)\qquad (u\to0).
\ee
In the appendix it is shown that (see (\ref{a2}))
\bee\label{e24}
x^2=\log\,s\bl(1+\frac{A(u)}{L}+\frac{B(u)}{2L^2}+\frac{C(u)}{4L^3}+O(L^{-4},(sL)^{-1})\br),
\ee
where we have defined
\[L=\log\,s,\quad L_1=\log\sqrt{\log\,s},\quad a=2\sqrt{\pi},\]
and the first few coefficients are
\[A(u)=-\log\,au-L_1,\quad B(u)=-A(u)-1,\quad C(u)=A^3(u)+3A(u)+\f{7}{2}.\]
We define the integrals
\[\la_n:=\int_0^\infty (\log\,u)^n e^{-u} du,\qquad (n=0, 1, 2, \ldots),\]
so that $\la_0=1$ and
\[\la_1=-\gamma,\quad \la_2=\gamma^2+\frac{\pi^2}{6},\quad \la_3=-\gamma^2-\frac{\pi^2\gamma}{2}-2\zeta(3),\]
where $\gamma=0.57721...$ is the Euler-Mascheroni constant and $\zeta$ is the Riemann zeta function.

Then, in the case $n=2$, we have upon extending the upper limit of integration to $\infty$ (thereby introducing an exponentially small error of O$(e^{-\al s})$) 
\begin{eqnarray}
I_{2,m}&\sim& \frac{\sqrt{\pi} \log\,s}{s} \int_0^\infty \bl(1+\frac{A(u)}{L}+\frac{B(u)}{2L^2}+\frac{C(u)}{4L^3}+\cdots\br) e^{-u}du\nonumber\\
&=&\frac{\sqrt{\pi} \log\,s}{s} \bl\{1+\frac{{\cal A}_2}{L}+\frac{{\cal B}_2}{2L^2}+\frac{{\cal C}_2}{4L^3}+ \cdots\br\}\qquad (m\to\infty),\label{e25}
\end{eqnarray}
where
\[{\cal A}_2=G-L_1,\qquad {\cal B}_2=-G+L_1-1,\]
\[{\cal C}_2=G^2+(3-2L_1)G+L_1^2-3L_1+\frac{\pi^2}{6}+\frac{7}{2}\]
with $G:=\gamma-\log\,a$. 

In Tricomi's case with $n=1$, we find from (\ref{e24}) that
\[x=\sqrt{\log\,s}\bl(1+\frac{A(u)}{2L}+\frac{2B(u)-A^2(u)}{8L^2}+\frac{2C(u)-2A(u)B(u)+A^3(u)}{16L^3}+O(L^{-4},(sL)^{-1})\br).\]
Hence we obtain the expansion
\bee\label{e26}
I_{1,m}\sim \frac{\sqrt{\pi \log\,s}}{s}\bl\{1+\frac{{\cal A}_1}{2L}-\frac{{\cal B}_1}{8L^2}+\frac{{\cal C}_1}{16L^3}+ \cdots\br\}\qquad (m\to\infty),
\ee
where
\[{\cal A}_1=G-L_1,\qquad {\cal B}_1=G^2+2(1-L_1)G+L_1^2-2L_1+\frac{\pi^2}{6}+2,\]
\[{\cal C}_1=G^3+(4-3L_1)G^2+(3L_1^2-8L_1+8+\frac{\pi^2}{6})G+4L_1^2-8L_1-L_1^3\]
\[+\frac{\pi^2}{6}(4-3L_1)+2\zeta(3)+7\]
and we recall that $s=m+1$. It is seen that the leading term in this expansion agrees with the result stated in (\ref{e10}). 

As a numerical verification of the expansions (\ref{e25}) and (\ref{e26}) we present\footnote{In the tables we write $x\times 10^y$ as $x(y)$.} in Table 1 the asymptotic values of $I_{n,m}$ ($n=1, 2$) for different $m$ compared with the values obtained by high-precision evaluation of (\ref{e11}).
In Table 2 we show the absolute relative error in the computation of $I_{1,m}$ for different $m$ as a function of the truncation index $k$ in the expansion (\ref{e26}).
\begin{table}[h]
\caption{\footnotesize{The values of $I_{n,m}$ ($n=1, 2$) compared with the expansions (\ref{e25}) and (\ref{e26}) for different $m$.}}
\begin{center}
\begin{tabular}{|c|ll|ll|}
\hline
&&&&\\[-0.3cm]
\mcol{1}{|c|}{$m$} & \mcol{1}{c}{$I_{1,m}$} & \mcol{1}{c|}{(2.7)} &  \mcol{1}{c}{$I_{2,m}$} & \mcol{1}{c|}{(2.6)}\\
[.1cm]\hline
&&&&\\[-0.25cm]
$10^2$ & $3.116097(-2)$ & $3.119672(-2)$ & $5.694564(-2)$ & $5.695077(-2)$\\
$10^3$ & $4.058838(-3)$ & $4.060226(-3)$ & $9.413132(-3)$ & $9.414250(-3)$\\
$10^4$ & $4.826833(-4)$ & $4.827422(-4)$ & $1.322800(-3)$ & $1.32285?(-3)$\\
$10^5$ & $5.494877(-5)$ & $5.495164(-5)$ & $1.710063(-4)$ & $1.710084(-4)$\\
$10^6$ & $6.094732(-6)$ & $6.094889(-6)$ & $2.101178(-5)$ & $2.101184(-5)$\\
[.1cm]\hline
\end{tabular}
\end{center}
\end{table}
\begin{table}[h]
\caption{\footnotesize{The absolute relative error in $I_{1,m}$ using the expansion (\ref{e26}) for different truncation index $k$.}}
\begin{center}
\begin{tabular}{|c|l|l|l|}
\hline
&&&\\[-0.3cm]
\mcol{1}{|c|}{$k$} & \mcol{1}{c|}{$m=10^4$} & \mcol{1}{c|}{$m=10^5$} &  \mcol{1}{c|}{$m=10^6$}\\
[.1cm]\hline
&&&\\[-0.25cm]
0 & $1.143(-1)$ & $9.447(-2)$ & $8.094(-2)$ \\
1 & $5.526(-3)$ & $3.686(-3)$ & $2.656(-3)$ \\
2 & $1.361(-5)$ & $1.013(-4)$ & $7.440(-5)$ \\
3 & $1.220(-5)$ & $5.214(-5)$ & $2.581(-5)$ \\
[.1cm]\hline
\end{tabular}
\end{center}
\end{table}

\vspace{0.6cm}

\begin{center}
{\bf 3.\ An extension of an integral arising in airfoil theory}
\end{center}
\setcounter{section}{3}
\setcounter{equation}{0}
\renewcommand{\theequation}{\arabic{section}.\arabic{equation}}
In airfoil theory, the expression for the downward velocity $w(a)$ at a point $x=a$ over the wingspan induced by a single vortex filament is \cite[p.~201]{PT}
\[w(a)=-\frac{1}{4\pi}\int_{-1}^1 \frac{1}{x-a}\,\frac{d\Gamma}{dx}\,dx\qquad(-1<a<1),\]
where $\Gamma$ is the circulation and the $x$-coordinate is normalised to the wingspan so that $-1\leq x\leq 1$. For an elliptical distribution $\Gamma=\Gamma_0 (1-x^2)^{1/2}$, where $\Gamma_0$ is the circulation in the middle of the airfoil, we have
\[w(a)=\frac{\Gamma_0}{4\pi}\int_{-1}^1 \frac{x}{x-a}\,\frac{dx}{\sqrt{1-x^2}}=\frac{\Gamma_0}{4\pi}\bl\{\int_{-1}^1\frac{dx}{\sqrt{1-x^2}}+a\int_{-1}^1 \frac{1}{x-a}\,\frac{dx}{\sqrt{1-x^2}}\br\}\]
\bee\label{e71}
=\frac{\Gamma_0}{4\pi}\bl\{\pi+a\int_{-1}^1 \frac{1}{x-a}\,\frac{dx}{\sqrt{1-x^2}}\br\}.
\ee
It is clear that $w(-a)=w(a)$ so that it is sufficient to consider $0<a<1$. 

The integral appearing in (\ref{e71}) must be treated as a Cauchy principal value; that is, we evaluate the integral up to $\epsilon$ either side of $x=a$ and then let $\epsilon\to 0$. Thus we have, with $R(w):=(1-a^2-2aw-w^2)^{1/2}$,
\[\int_{-1}^1 \frac{1}{x-a}\,\frac{dx}{\sqrt{1-x^2}}=\int_\epsilon^{1-a}\frac{dw}{w R(w)}+\int_{1+a}^\epsilon \frac{dw}{w R(-w)}\]
\[=\frac{1}{\sqrt{1-a^2}}\bl\{\bl[\log\,\frac{w}{2(1\!-\!a^2\!-\!aw\!+\!\sqrt{1-a^2}R(w))}\br]^{1-a}_\epsilon+
\bl[\log\,\frac{w}{2(1\!-\!a^2\!+\!aw\!+\!\sqrt{1-a^2}R(-w))}\br]^\epsilon_{1+a}\br\}\]
\[=-\frac{2a\epsilon}{(1-a^2)^{3/2}}+O(\epsilon^2)\]
after some routine algebra. Thus the integral vanishes as $\epsilon\to0$ with the result that 
\bee\label{e71b}
\int_{-1}^1 \frac{x}{x-a}\,\frac{dx}{\sqrt{1-x^2}}=\pi
\ee
independent of $a$ \cite[p.~202]{PT}. Hence $w(a)=\Gamma_0/4$, which
means that in the case of an elliptical distribution of the circulation the induced downward velocity is constant across the wing.

The situation when $d\Gamma/dx=x^{2n+1} \Gamma_0/\sqrt{1-x^2}$, $n=0, 1, 2, \ldots\ $, corresponds to a flattening of the basic elliptic profile ($n=0$). For example, when $n=1, 2$ we have the profiles
\[\Gamma=\frac{\Gamma_0}{2} (2+x^2) \sqrt{1-x^2}\ \ \ (n=1), \qquad \Gamma=\frac{\Gamma_0}{15} (8+4x^2+3x^4) \sqrt{1-x^2}\ \ \ (n=2);\]
as $n$ increases the profile becomes progressively flatter in the central portion of the wing. Then the extension of the integral (\ref{e71}) we consider is
\bee\label{e72} 
J_n(a;\mu):=\int_{-1}^1 \frac{x^{2n+1}}{x-a}\,\frac{dx}{(1-x^2)^\mu}\qquad (\mu<1,\ n=0, 1, 2, \ldots).
\ee
which reduces to (\ref{e71}) when $n=0$ and $\mu=\fs$. It is easily seen that $J_n(-a;\mu)=J_n(a;\mu)$, so that again it is sufficient to consider $0<a<1$.
\vspace{0.5cm}

\noindent{\bf 3.1\ \ The evaluation of the integral $J_n(a;\mu)$}
\vspace{0.2cm}

\noindent
We write $J_n(a;\mu)$ in the form
\[J_n(a;\mu)=\int_{-1}^1\frac{x^{2n+1}-a^{2n+1}}{x-a}\,\frac{dx}{(1-x^2)^\mu}+a^{2n+1}\int_{-1}^1\frac{1}{x-a}\,\frac{dx}{(1-x^2)^\mu}\equiv I_1+a^{2n+1}I_2.\]
In the first integral we employ the expansion
\bee\label{e71a}
\frac{x^{m+1}-a^{m+1}}{x-a}=\sum_{r=0}^{m}a^{m-r}x^r\qquad (m=0, 1, 2, \ldots)
\ee
to find
\[I_1=\sum_{r=0}^n a^{2n-2r}\int_{-1}^1\frac{x^{2r}}{(1-x^2)^\mu}\,dx=
\frac{a^{2n}\sqrt{\pi} \g(1-\mu)}{\g(\f{3}{2}-\mu)} 
\sum_{r=0}^n \frac{(\fs)_r a^{-2r}}{(\f{3}{2}-\mu)_r}\]
\[=\frac{a^{2n}\sqrt{\pi} \g(1-\mu)}{\g(\f{3}{2}-\mu)}\,{}_2F_1(\fs,1;\f{3}{2}-\mu;a^{-2})|_n,\] 
where the subscript $n$ on the Gauss hypergeometric function ${}_2F_1$ denotes that only the first $n+1$ terms are to be taken.

The second integral is
\[I_2=\sum_{k=0}^\infty\frac{(\mu)_k}{k!} \int_{-1}^1\frac{x^{2k}}{x-a}\,dx
=\sum_{k=0}^\infty\frac{(\mu)_k}{k!}\bl\{a^{2k}\int_{-1}^1\frac{dx}{x-a}+\int_{-1}^1\frac{x^{2k}-a^{2k}}{x-a}\,dx\br\}\]
\[=\frac{1}{(1-a^2)^\mu} \log \bl(\frac{1-a}{1+a}\br)+\sum_{k=1}^\infty\frac{(\mu)_k}{k!}\int_{-1}^1\frac{x^{2k}-a^{2k}}{x-a}\,dx.\]
Use of (\ref{e71a}) in the above integral shows that
\[S:=\sum_{k=1}^\infty\frac{(\mu)_k}{k!}\int_{-1}^1\frac{x^{2k}-a^{2k}}{x-a}\,dx=\frac{1}{a}\sum_{k=1}^\infty\frac{(\mu)_k}{k!} \sum_{r=0}^{k-1} a^{2k-2r}\int_{-1}^1 x^{2r}dx\]
\[=2a\sum_{k=0}^\infty\frac{(\mu)_{k+1}a^{2k}}{(k+1)!}\sum_{r=0}^k\frac{a^{-2r}}{2r+1}=
2a\sum_{r=0}^\infty \frac{a^{-2r}}{2r+1} \sum_{k=r}^\infty \frac{(\mu)_{k+1}a^{2k}}{(k+1)!}\]
\[=2a\sum_{r=0}^\infty\frac{1}{2r+1}\sum_{k=0}^\infty\frac{(\mu)_{k+r+1}a^{2k}}{(k+r+1)!},\]
where we have put $k\to k+r$ in the last sum. From the identity $(\mu)_{k+r+1}=(\mu)_{r+1} (\mu+r+1)_k$, we then obtain
\begin{eqnarray*}
S&=&2a\sum_{r=0}^\infty \frac{(\mu)_{r+1}}{(2r+1) (r+1)!}\sum_{k=0}^\infty \frac{(\mu+r+1)_ka^{2k}}{(r+2)_k}\\
&=&2a\sum_{r=1}^\infty \frac{(\mu)_{r}}{(2r-1) r!}\,{}_2F_1(1,\mu+r;r+1;a^2)\\
&=&\frac{2a}{1-a^2} \sum_{r=1}^\infty \frac{(\mu)_{r}}{(2r-1) r!}\,{}_2F_1\bl(1,1-\mu;r+1;-\frac{a^2}{1-a^2}\br)
\end{eqnarray*}
upon use of the well-known Euler transformation \cite[(15.8.1)]{DLMF} for the hypergeometric function.

Then we have the result
\[
J_n(a;\mu)=\frac{a^{2n}\sqrt{\pi} \g(1-\mu)}{\g(\f{3}{2}-\mu)}\,{}_2F_1(\fs,1;\f{3}{2}-\mu;a^{-2})|_n
+\frac{a^{2n+1}}{(1-a^2)^\mu}\,\log \bl(\frac{1-a}{1+a}\br)\]
\bee\label{e74}
+\frac{2a^{2n+2}}{1-a^2} \sum_{r=1}^\infty \frac{(\mu)_{r}}{(2r-1) r!}\,{}_2F_1\bl(1,1-\mu;r+1;-\frac{a^2}{1-a^2}\br)
\ee
valid for $\mu<1$ and non-negative integer $n$. We remark that the convergence of the infinite sum in (\ref{e74}) is slow, being controlled by $r^{\mu-2}$ as $r\to\infty$. A simple means of accelerating the convergence is given in the sub-section below.
\vspace{0.5cm}

\noindent{\bf 3.2\ \ Improved convergence of (\ref{e74})}
\vspace{0.2cm}

\noindent
We can write the hypergeometric function appearing in the infinite sum in (\ref{e74}) in the form
\[{}_2F_1(1,1-\mu;r+1;-X)=1-\frac{(1-\mu)X}{r+1}+\frac{(1-\mu)_2X^2}{(r+1)_2}{}_2F_1(1,3-\mu;r+3;-X),\]
where, for convenience, we have put $X:=a^2/(1-a^2)$. Then
\[S=\frac{2a}{1-a^2}\bl(\sigma_0-(1-\mu)X \sigma_1+(1-\mu)_2X^2 \sigma_2\br)\]
\bee\label{e75}
+\frac{2a}{1-a^2} (1-\mu)_2X^2 \sum_{r=1}^\infty \frac{(\mu)_r}{(2r-1) (r+2)!}\,\bl\{{}_2F_1(1,3-\mu;r+3;-X)-1\br\},
\ee
where
\begin{eqnarray*}
\sigma_0&=&\sum_{r=1}^\infty\frac{(\mu)_r}{(2r-1) r!}=1-\frac{\sqrt{\pi}\g(1-\mu)}{\g(\fs-\mu)},\\
\sigma_1&=&\sum_{r=1}^\infty\frac{(\mu)_r}{(2r-1) (r+1)!}=\frac{2-3\mu}{3(1-\mu)}-\frac{2\sqrt{\pi}\g(1-\mu)}{3\g(\fs-\mu)},\\
\sigma_2&=&\sum_{r=1}^\infty\frac{(\mu)_r}{(2r-1) (r+2)!}=\frac{16-25\mu+15\mu^2}{30(1-\mu)(2-\mu)}-
\frac{4\sqrt{\pi}\g(1-\mu)}{15\g(\fs-\mu)}.
\end{eqnarray*}

The convergence of the final sum in (\ref{e75}) is now controlled by terms of O$(r^{\mu-5})$ as $r\to\infty$, which is an improvement 
on that in (\ref{e74}). This procedure can be continued to produce higher rates of convergence, but at the expense of the evaluation of higher-order sums $\sigma_m$, $m\geq 2$.

\vspace{0.5cm}

\noindent{\bf 3.3\ \ The special case $n=0$, $\mu=\fs$}
\vspace{0.2cm}

\noindent
We demonstrate that the result in (\ref{e74}) reduces to the value given in (\ref{e71b}) when $n=0$ and $\mu=\fs$.
We have from (\ref{e74}) with $X=a^2/(1-a^2)$
\[J_0(a;\fs)=\pi+\frac{a}{\sqrt{1-a^2}}\,\log \bl(\frac{1-a}{1+a}\br)+2X\sum_{r=1}^\infty \frac{(\mu)_r}{(2r-1) r!}\,{}_2F_1(\fs,1;r+1;-X). \]
The sum can be written in the form
\[T:=X\sum_{r=0}^\infty\frac{(\fs)_{r+1}}{(r+\fs) (r+1)!} \sum_{k=0}^\infty\frac{(\fs)_k (-X)^k}{(r+2)_k}
=X\sum_{k=0}^\infty (\fs)_k (-X)^k \sum_{r=0}^\infty \frac{(\fs)_{r+1}}{(r+\fs) (r+k+1)!}\]
\[=X\sum_{k=0}^\infty \frac{(\fs)_k (-X)^k}{(k+1)!}\,{}_2F_1(\fs,1;k+2;1).\]
Application of Gauss' summation theorem  series \cite[(15.4.20)]{DLMF} shows that the ${}_2F_1(1)$ series has the value
$(k+1) \g(k+\fs)/\g(k+\f{3}{2})$, whence
\[T=2X\sum_{k=0}^\infty \frac{(\fs)_k (\fs)_k}{(\f{3}{2})_k k!}\,(-X)^k=2X\,{}_2F_1(\fs,\fs;\f{3}{2};-X)\]
\[=\frac{2X}{\sqrt{1+X}}\,{}_2F_1\bl(\fs,1;\f{3}{2};\frac{X}{1+X}\br)=\frac{a}{\sqrt{1-a^2}}\,\log \bl(\frac{1+a}{1-a}\br)\]
upon use of Euler's transformation followed by evaluation of the resulting hypergeometric function by \cite[(15.4.2)]{DLMF}.
Hence, we recover the result
\[J_0(a;\fs)=\pi\]
as given in (\ref{e71b}).

\vspace{0.6cm}

\begin{center}
{\bf Appendix: \ The inversion of (\ref{e23inv})}
\end{center}
\setcounter{section}{1}
\setcounter{equation}{0}
\renewcommand{\theequation}{\Alph{section}.\arabic{equation}}
We require the inversion of the expression
\[\mbox{erfc}\,x=\frac{2u}{s}\bl(1-\frac{u}{2s}+\cdots\br)\qquad (u\to0).\]
Since the limit $u\to0$ corresponds to $x\to+\infty$, we employ the asymptotic expansion \cite[(7.12.1)]{DLMF}
\[\mbox{erfc}\,x=\frac{e^{-x^2}}{\sqrt{\pi} x} \bl\{1-\frac{1}{2x^2}+\frac{3}{4x^4}+O(x^{-6})\br\}\qquad (x\to+\infty)\]
to yield
\[xe^{x^2}=\frac{s(1-\fs x^{-2}+\f{3}{4}x^{-4}-\cdots)}{au(1-u/(2s)+\cdots)},\qquad a=2\sqrt{\pi}.\]
Taking logarithms, we obtain
\bee\label{a1}
x^2=\log\,s \bl(1-\frac{\log\,au}{L}-\frac{\log\,x}{L}-\frac{1}{2x^2L}+\frac{5}{8x^4L}+\cdots +O((sL)^{-1})\br),
\ee
where we set $L=\log\,s$, $L_1=\log \sqrt{\log\,s}$.

Following the iterative procedure described in \cite[pp.~25--26]{NB} we have as a first approximation $x^2=\log\,s$, whence
\[x^2=\log\,s \bl(1+\frac{A(u)}{L}\br),\qquad A(u)=-\log\,au-L_1.\]
Then with
\[\frac{\log\,x}{L}=\frac{L_1}{L}+\frac{A(u)}{2L^2}+\ldots ,\]
we find a second approximation given by
\[x^2=\log\,s \bl(1+\frac{A(u)}{L}+\frac{B(u)}{2L^2}+\cdots \br),\qquad B(u)=-A(u)-1.\]
This last result produces
\[\frac{\log\,x}{L}=\frac{L_1}{L}+\frac{A(u)}{2L^2}+\frac{B(u)-A^2(u)}{4L^3}+\cdots ,\]
so that we obtain the third approximation 
\bee\label{a2}
x^2=\log\,s \bl(1+\frac{A(u)}{L}+\frac{B(u)}{2L^2}+\frac{C(u)}{4L^3}+O(L^{-4}, (sL)^{-1})\br),
\ee
where
\[C(u)=A^2(u)+3A(u)+\f{7}{2}.\] 
\vspace{0.6cm}

\noindent{ {\bf Acknowledgement:} \, The author wishes to acknowledge F. Mainardi for bringing to his attention the 1933 paper by Tricomi.

\end{document}